\documentclass[11pt]{article}
\usepackage{amsmath,amsfonts}
\usepackage{verbatim}
\usepackage{latexsym}
\usepackage{graphicx}
\usepackage{dsfont}
\usepackage{indentfirst}
\usepackage{xcolor}
\usepackage{booktabs}
\makeatletter
\@addtoreset{equation}{section}

\begin{document}

\newcommand{\E}{\mathbb{E}}
\newcommand{\PP}{\mathbb{P}}
\newcommand{\RR}{\mathbb{R}}
\newcommand{\LL}{\mathbb{L}}

\newtheorem{theorem}{Theorem}[section]
\newtheorem{lemma}[theorem]{Lemma}
\newtheorem{coro}[theorem]{Corollary}
\newtheorem{defn}[theorem]{Definition}
\newtheorem{assp}[theorem]{Assumption}
\newtheorem{expl}[theorem]{Example}
\newtheorem{prop}[theorem]{Proposition}
\newtheorem{rmk}[theorem]{Remark}

\newcommand\tq{{\scriptstyle{3\over 4 }\scriptstyle}}
\newcommand\qua{{\scriptstyle{1\over 4 }\scriptstyle}}
\newcommand\hf{{\textstyle{1\over 2 }\displaystyle}}
\newcommand\hhf{{\scriptstyle{1\over 2 }\scriptstyle}}

\newcommand{\proof}{\noindent {\it Proof}. }
\newcommand{\eproof}{\hfill $\Box$} 

\def\a{\alpha} \def\g{\gamma}
\def\e{\varepsilon} \def\z{\zeta} \def\y{\eta} \def\ot{\theta}
\def\vo{\vartheta} \def\k{\kappa} \def\lbd{\lambda} \def\m{\mu} \def\n{\nu}
\def\x{\xi}  \def\r{\rho} \def\s{\sigma}
\def\p{\phi} \def\f{\varphi}   \def\w{\omega}
\def\q{\surd} \def\i{\bot} \def\h{\forall} \def\j{\emptyset}

\def\be{\beta} \def\de{\delta} \def\up{\upsilon} \def\eq{\equiv}
\def\ve{\vee} \def\we{\wedge}

\def\F{{\cal F}}
\def\T{\tau} \def\G{\Gamma}  \def\D{\Delta} \def\O{\Theta} \def\L{\Lambda}
\def\X{\Xi} \def\S{\Sigma} \def\W{\Omega}
\def\M{\partial} \def\N{\nabla} \def\Ex{\exists} \def\K{\times}
\def\V{\bigvee} \def\U{\bigwedge}

\def\1{\oslash} \def\2{\oplus} \def\3{\otimes} \def\4{\ominus}
\def\5{\circ} \def\6{\odot} \def\7{\backslash} \def\8{\infty}
\def\9{\bigcap} \def\0{\bigcup} \def\+{\pm} \def\-{\mp}
\def\la{\langle} \def\ra{\rangle}

\def\tl{\tilde}
\def\trace{\hbox{\rm trace}}
\def\diag{\hbox{\rm diag}}
\def\for{\quad\hbox{for }}
\def\refer{\hangindent=0.3in\hangafter=1}

\newcommand\wD{\widehat{\D}}

\title{
\bf Delay-dependent Asymptotic Stability of Highly Nonlinear Stochastic Differential Delay Equations Driven by $G$-Brownian Motion
 }

\author{
{\bf
Chen Fei${}^{2}$,~
Weiyin Fei${}^{1}$\thanks{Corresponding author. E-mail:  wyfei@ahpu.edu.cn},~
Xuerong Mao${}^{3}$,~
Litan Yan${}^{2}$
 }
\\
${}^1$ School of Mathematics and Physics, \\
Anhui Polytechnic University, Wuhu, Anhui 241000, China.\\
${}^2$ Department of Statistics, College of Science,\\
Donghua University, Shanghai, 200051, China.\\
${}^3$Department of Mathematics and Statistics,\\
 University of Strathclyde, Glasgow G1 1XH, UK.\\
}
\date{}

\maketitle

\begin{abstract}
     Based on the classical probability, the stability of stochastic differential delay equations (SDDEs) whose coefficients are either linear or nonlinear but bounded by linear functions have been investigated intensively. Moreover, the delay-dependent stability of highly nonlinear hybrid stochastic differential equations has also been studied recently. In this paper, by using the nonlinear expectation theory, we explore the delay-dependent stability of a class of highly nonlinear hybrid stochastic differential delay equations driven by $G$-Brownian motion ($G$-SDDEs). Firstly, we give preliminaries on sublinear expectation. Then, the delay-dependent criteria of the stability and boundedness of solutions to $G$-SDDEs is developed. Finally, an illustrative example is analyzed by the $\varphi$-max-mean algorithm.


\medskip \noindent{\small \bf MR Subject Classification:} 60H10.
\end{abstract}

\section{Introduction}

A sublinear expectation (c.f., \cite{P2006,P2008,P2010}) can be represented as the upper expectation of
a subset of linear expectations. In
most cases, this subset is often treated as an uncertain model of probabilities. Peng introduced the $G$-expectation theory ($G$-framework) (see \cite{P2006} and the references
therein) in 2006, where the notion of
$G$-Brownian motion and the corresponding stochastic calculus of It\^o's type were established.  Moreover, the stochastic calculus on $G$-Brownian motion is explored by Denis {\sl et al.} \cite{DHP}, Fei and Fei \cite{FF}, Hu and Peng \cite{HP}, Li and Peng \cite{LP}, Peng and Zhang \cite{PZ2017}, Soner et al. \cite{STZ2} and Zhang \cite{Z2016}, etc.

So far, there is large amount of literature on the problem of asset pricing and financial decisions under model uncertainty. Chen and Epstein \cite{CE} put forward to the model of an intertemporal recursive utility, where risk and ambiguity are differentiated, but uncertainty is only a mean uncertainty without a volatility uncertainty. The model of the optimal consumption and portfolio with ambiguity are also investigated in Fei \cite{FeiINS,FeiSM}. Epstein and Ji \cite{EJ1,EJ2} generalized the Chen-Epstein model and maintained a separation between risk aversion
and intertemporal substitution. We know that equivalence of priors is an optional assumption in Gilboa and Schmeidler \cite{GS}. Apart from very recent developments, the stochastic calculus presumes a probability space framework. However, from an economics perspective, the assumption of equivalence seems far from innocuous.
Informally, if her environment is complex, how could the decision-maker come to
be certain of which scenarios regarding future asset prices and rates of return, for
example, are possible? In particular, ambiguity about volatility implies ambiguity
about which scenarios are possible, at least in a continuous time setting.
A large amount of literature has argued that the stochastic time varying volatility is important for understanding features of asset returns, and particularly empirical regularities in the derivative markets.

On the other hand, it is well to be known that the study of differential equations from both the theory and its applications is important. The classical stochastic differential equation with Brownian motion doesn't take an ambiguous factor into consideration. Thus, in some complex environments, these equations are too restrictive for describing some phenomena. Recently, under uncertainty, a kind of stochastic differential equations driven by $G$-Brownian motion has been investigated by Deng {\sl et al.} \cite{DFFM}, Fei and Fei \cite{FF2019}, Fei {\sl et al.} \cite{FFY}, Gao \cite{G}, Hu {\sl et al.} \cite{HJPS},  Li and Yang \cite{LY}, Lin\cite{Lq2013},  Lin \cite{Lin2013}, Luo and Wang \cite{LW2014}, Peng \cite{P2010}, Ren {\sl et al.} \cite{RYS}, Xu {\sl et al} \cite{XGH}, Yao and Zong \cite{YZ}, and Yin {\sl et al.} \cite{YCR}.

We know that the stability of the classical stochastic differential equations is an important topic in the study of stochastic systems. Since Lyapunov \cite{Ly} and LaSalle \cite{La} obtained the stability for nonlinear systems, with the help of It\^o's breakthrough work about It\^o calculus, Hasminskii \cite{Has} first studied the stability of the linear It\^o equation perturbed a noise. Since then, the stability analysis for stochastic differential systems has been done by many researchers.
 The stabilization and destabilization of hybrid systems of stochastic differential equations is explored by many researchers, such as Mao et al. \cite{MYY}. In Luo and Liu \cite{LL}, the stability of infinite dimensional stochastic evolution equations with memory and Markovian jumps is investigated.

Since a system described by the stochastic differential equation driven by $G$-Brownian motion provides a characterization of the real world with both randomness and ambiguity, it is necessary to investigate its stability like a classical stochastic differential equation. One kind of exponential stability for stochastic differential equations driven by $G$-Brownian motion is discussed by Zhang and Chen \cite{ZC} where the quasi-sure analysis is used.  Fei and Fei \cite {FF1} investigated the quasi sure exponential stability by $G$-Lyapunov functional method for obtaining different results.

However, in many real systems, such as science, industry, economics and finance etc., we will run into time lag.
Differential delay equations (DDEs) have been used to set up
such time-delay systems. However, since the time-delay often causes the instability of systems, the stability of DDEs  has  been  investigated intensively for more than 50 years.
In 1980's, the stochastic delay differential equations were developed in order to model real systems
which are subject to external noises.
Since then, in the study of SDDEs including hybrid SDDEs the stability analysis has been one of the most important topics (see, e.g., \cite{LM,Liu12,Mao2007}).

The existing delay-dependent stability criteria are mainly provided for the hybrid SDDEs whose coefficients are either linear or nonlinear but bounded by linear functions (or, satisfy the linear growth condition).  Recently, \cite{HMS,HMZ} initiate the investigation on the stability of the hybrid highly nonlinear stochastic delay differential equations driven by Brownian noise. Based on highly nonlinear hybrid SDDEs (see, e.g., \cite{HMS,HMZ}),
\cite{FHMS} has further established the delay-dependent stability criterion. Recently, the delay-dependent stability of the highly nonlinear hybrid  SDDEs is investigated in \cite{FFMXY,FSFMY, FHMS2,MFFM,SFFM2019,SFFM2020}.

 Moreover, to our best knowledge, the current states of the considered real systems mainly follow an SDDE disturbed by Brownian motion. In fact, due to uncertainty including probabilistic
  and Kightian ones \cite{CE}, many real systems may be disturbed by $G$-Brownian motion. Thus we can characterize our some real systems by Peng's sublinear expectation framework. Moreover, in this paper, we will explore the dependent stability criteria of SDDEs driven by $G$-Brownian motion ($G$-SDDEs), where the coefficients of $G$-SDDEs are highly nonlinear without linear growth conditions. Our results extend the ones under linear expectation framework discussed in Fei {\sl et al.}  \cite{FHMS} .

We now give an example on the stability of the following  highly nonlinear $G$-SDDE
\begin{align}\label {$G$-SDDE1}
dX(t) = [-X^3(t)-X(t-\de(t))]dt+\frac12e^{-t}dB(t),
\end{align}
where $\de(t)\geq0$ is a time lag, $X(t)\in \mathbb{R}$ is the state of the highly non-linear hybrid system, $B(t)$ is a scalar $G$-Brownian motion (related notions see \cite{P2010}).

For the above system (\ref{$G$-SDDE1}), if the time delay $\de(t)=0.01$, the computer simulation shows it is asymptotically stable (see Figure 4.1). If the time-delay is large, say $\de(t)=2$,  the computer simulation shows that $G$-SDDE (\ref{$G$-SDDE1}) is unstable (see Figure 4.2). In other words, whether the delay $G$-SDDE is stable or not depends on how small or large the time-delay is.  On the other hand, both drift and diffusion coefficients of $G$-SDDE affect the stability of the system due to highly nonlinear. However, there is no delay dependent criterion which can be applied to $G$-SDDE to derive a sufficient bound on the time-delay $\de(t)$ such that $G$-SDDE is stable. The aim of this paper is to establish the delay dependent criteria for highly nonlinear $G$-SDDEs.

The main contributions of this paper are presented as follows:

(i). We first try to give the criteria on delay-dependent stability of $G$-SDDEs with highly nonlinear coefficients.

(ii). The delay-dependent criteria are established by new mathematical techniques, such as constructing Lyapunov functional.

(iii). The stochastic calculus on $G$-Brownian motion is applied to solve the stability of the systems with ambiguity.

 The arrangement of the paper is as follows. In Section 2, we give preliminaries on sublinear expectations and $G$-Bwownian motions. Furthermore, we characterize
 the properties of $G$-Brownian motions and $G$-martingales.
In Section 3, the quasi sure exponential stability of the solutions to $G$-SDDE is studied. Next, we give an illustrative example in Section 4, where we use the $\varphi$-max-mean algorithm suggested by Peng \cite{P-scm}. Finally, Section 5 concludes this paper.

\section{Notations and Preliminaries}

 In this section, we first give the notion of sublinear expectation space $(\Omega, {\cal H}, \hat{\mathbb{E}})$, where $\Omega$ is a given state set and $\cal H$ a linear space of real valued functions defined on $\Omega$. The space $\cal H$ can be considered as the space of random variables. The following concepts come from  Peng \cite{P2010}.

\begin{defn} \label{Def2.1}
 A sublinear expectation $\hat{\mathbb E}$ is a functional $\hat{\mathbb  E}$: ${\cal H}\rightarrow {\mathbb  R}$ satisfying

\noindent(i) Monotonicity: $\hat{\mathbb{E}}[X]\geq \hat{\mathbb{E}}[Y]$ if $X\geq Y$;

\noindent(ii) Constant preserving: $\hat{\mathbb{E}}[c]=c$;

\noindent (iii) Sub-additivity: For each $X,Y\in {\cal H}$, $\hat{\mathbb{E}}[X+Y]\leq \hat{\mathbb{E}}[X]+\hat{\mathbb{E}}[Y]$;

\noindent(iv) Positivity homogeneity: $\hat{\mathbb{E}}[\lambda X]=\lambda \hat{\mathbb{E}}[X]$ for $\lambda\geq 0$.
\end{defn}

\begin{defn} \label{Def2.2} Let $(\Omega,{\cal H}, \hat{\mathbb{E}})$ be a sublinear expectation space. $(X(t))_{t\geq0}$ is called a $d$-dimensional stochastic process if for each $t\geq0$, $X(t)$ is a $d$-dimensional random vector in $\cal H$.

A $d$-dimensional process $(B(t))_{t\geq0}$ on a sublinear expectation space $(\Omega,{\cal H}, \hat{\mathbb{E}})$ is called a $G$-Brownian motion if the following properties are satisfies:

\noindent(i) $B_0(\omega)=0$;

\noindent (ii) for each $t,s\geq 0$, the increment $B(t+s)-B(t)$ is $N(\{0\}\times s\Sigma)$-distributed and is independent from $(B(t_1),B(t_2),\cdots, B(t_n))$, for each $n\in {\mathbb N}$ and $0\leq t_1\leq \cdots\leq t_n\leq t$.

\end{defn}

More details on the notions of $G$-expectation ${\hat{\mathbb E}}$ and $G$-Brownian motion on the sublinear expectation space $(\Omega,{\cal H}, {\hat{\mathbb{E}}})$ may be found in Peng \cite{P2010}.

We now give the definition of It\^o integral. For simplicity, in the rest of the paper, we introduce It\^o integral with respect to 1-dimensional $G$-Brownian motion with $G(\alpha):=\frac{1}{2}\hat{\mathbb{E}}[\alpha B(1)^2]=\frac{1}{2}(\bar{\sigma}^2\alpha^+-\underline{\sigma}^2\alpha^-)$, where $\hat{\mathbb{E}}[B(1)^2]=\bar{\sigma}^2, {\cal E}[B(1)^2]:=-\hat{\mathbb{E}}[- B(1)^2]=\underline{\sigma}^2$,
$0<\underline{\sigma}\leq \bar{\sigma}<\infty$.

Let $p\geq1$ be fixed. We consider the following type of simple processes: for a given partition $\pi_T=(t_0,\cdots,t_N)$
of $[0, T]$, where $T$ can take $\infty$, we get
$$\eta_t(\omega)=\sum_{k=0}^{N-1}\xi_k(\omega)I_{[t_k,t_{k+1})}(t),$$
where $\xi_k\in L^p_G(\Omega_{t_k}), k=0,1,\cdots,N-1$ are given. The collection of these processes is denoted by $M_G^{p,0}(0,T)$. We denote by $M_G^p(0,T)$ the completion of $M_G^{p,0}(0,T)$  with the norm
 $$\|\eta\|_{M_G^p(0,T)}:=\left\{\int_0^T\hat{\mathbb{E}}|\eta(t)|^pdt\right\}^{1/p}<\infty.$$
 We provide the following property which comes from Denis et al. \cite{DHP} or Zhang and Chen \cite{ZC}.

 \vskip12pt

\begin{prop}\label{Pro1} Let $\hat{\mathbb{E}}$ be $G$-expectation. Then there exists a weakly compact family of probability measures $\cal P$ on $(\Omega, {\cal B}(\Omega))$ such that for all $X\in {\cal H}, \hat{\mathbb{E}}[X]=\max_{P\in{\cal P}}E^P[X]$, where $E^P[\cdot]$ is the linear expectation with respect to $P$.
\end{prop}

From the above Proposition \ref{Pro1}, we know that the weakly compact family of probability measures $\cal P$ characterizes the degree of Knightian uncertainty. Especially, if $\cal P$ is singleton, i.e. $\{P\}$, then the model has no ambiguity. The related calculus reduces to a classical one.
 
   The definition of stochastic integral $\int_0^T\eta(t) dB(t)$ can see Peng \cite{P2010}.
By \cite[Proposition 2.5]{FF} and  Peng \cite[Lemma 3.5 on page 43]{P2010}, we easily obtain the following lemma.

\begin{lemma}\label{BBI} For all $s< t\in[0,T]$, we have
  \begin{align*}
 & \hat{\mathbb E}\Big|\int_s^t\zeta(v)d<B(v)>\Big|^2\leq(t-s)\bar\sigma^2\hat{\mathbb E}\Big(\int_s^t|\zeta(v)|^2dv\Big),\\
 &\hat{\mathbb E}\Big|\int_s^t\zeta(v)dB(v)\Big|^2\leq \bar\sigma^2\hat{\mathbb E}\Big(\int_s^t|\zeta(v)|^2dv\Big).\\
\end{align*}
\end{lemma}

The following Burkholder-Davis-Gundy inequality provides the explicit bounds based on those in Gao \cite[Theorems 2.1-2.2]{G}.

\begin{lemma}(Burkholder-Davis-Gundy inequality)\label{BDG} Let $p>0$ and $\zeta=\{\zeta(s),s\in[0,T]\}\in M_G^p(0,T)$. Then, for all $t\in[0,T]$,
  \begin{align*}
 & \hat{\mathbb E}\sup\limits_{s\leq u\leq t}\Big|\int_s^u\zeta(v)d<B(v)>\Big|^p\\
 &\quad\leq(t-s)^{p-1}C_1(p,\bar\sigma)\hat{\mathbb E}\Big(\int_s^t|\zeta(v)|^2dv\Big)^{p/2},\\
 &\hat{\mathbb E}\sup\limits_{s\leq u\leq t}\Big|\int_s^u\zeta(v)dB(v)\Big|^p\leq C_2(p,\bar\sigma)\hat{\mathbb E}\Big(\int_s^t|\zeta(v)|^2dv\Big)^{p/2},\\
\end{align*}
where $C_1(p,\bar{\sigma})=\bar{\sigma}^p, C_2(p,\bar{\sigma})=C_p\bar{\sigma}^p$, and the constant $C_p$ is defined as follows
 \begin{align*}
 C_p&=\Big(\frac{32}{p}\Big)^{p/2},\quad\mbox{\rm if}~0<p<2;\\
 C_p&=4,\quad\mbox{\rm if}~p=2;\\
 C_p&=\Big(\frac{p^{p+1}}{2(p-1)^{p-1}}\Big)^{p/2},\quad\mbox{\rm if}~p>2.
\end{align*}

\end{lemma}

\textbf{Proof}: Based on the idea from Peng {\sl et al.} \cite{PSZ}, we can characterize the probability measure family $\mathcal P$ as follows:
$${\mathcal P}=\{P^{(\sigma_\cdot)}; (\sigma_t)_{t\in[0,\infty)}~\mbox{is a continuous process}, \sigma_t\in [\underline{\sigma},\bar\sigma]\}.$$
For each probability measure $P^{(\sigma_\cdot)}\in {\mathcal P}$, we have $dB(t)=\sigma_tdw(t)$, where $(w(t))_{t\ge0}$ is a standard Brownian motion under $P^{(\sigma_\cdot)}$. The corresponding linear expectation is denoted by $ E^{P^{\sigma_\cdot}}$. Thus, by the H\"old inequality, we get
 \begin{align*}
 & \hat{\mathbb E}\sup\limits_{s\leq u\leq t}\Big|\int_s^u\zeta(v)d<B(v)>\Big|^p\\
 & =\sup\limits_{P^{\sigma_\cdot}\in{\cal P}}E^{P^{\sigma_\cdot}}\sup\limits_{s\leq u\leq t}\Big|\int_s^u\zeta(v)\sigma_vd<w(v)>\Big|^p\\
 &\quad\leq(t-s)^{p-1}\bar{\sigma}^p\sup\limits_{P^{\sigma_\cdot}}E^{P^{\sigma_\cdot}}\Big(\int_s^t|\zeta(v)|^2dv\Big)^{p/2}\\
 &\quad=(t-s)^{p-1}\bar{\sigma}^p\hat{\mathbb E}\Big(\int_s^t|\zeta(v)\sigma_v|^2dv\Big)^{p/2}.
\end{align*}
From the classical Burkholder-Davis-Gundy inequality (see, e.g. Mao and Yuan \cite[Theorem 2.13 on page 70]{MY2006}), we get
\begin{align*}
 &\hat{\mathbb E}\sup\limits_{s\leq u\leq t}\Big|\int_s^u\zeta(v)dB(v)\Big|^p\\
 &=\sup\limits_{P^{\sigma_\cdot}\in{\cal P}}E^{P^{\sigma_\cdot}}\sup\limits_{s\leq u\leq t}\Big|\int_s^u\zeta(v)\sigma_vdw(v)\Big|^p\\
 &\leq C_p\sup\limits_{P^{\sigma_\cdot}\in{\cal P}}E^{P^{\sigma_\cdot}}\sup\limits_{s\leq u\leq t}\Big|\int_s^u|\zeta(v)\sigma_v|^2dv\Big|^{p/2}\\
& \leq C_p\bar\sigma^p\hat{\mathbb E}\Big(\int_s^t|\zeta(v)|^2dv\Big)^{p/2}.
\end{align*}
Thus, the proof is complete. $\Box$

 We now define $G$-upper capacity $\mathbb{V}(\cdot)$ and $G$-lower capacity $\mathcal{V}(\cdot)$ by
\begin{align*}
{\mathbb V}(A)=&\sup_{P\in {\cal P}}P(A),~\forall~ A\in {\cal B}(\Omega),\\
{\mathcal V}(A)=&\inf_{P\in {\cal P}}P(A),~\forall~ A\in {\cal B}(\Omega).
\end{align*}
Thus a property is called to hold quasi surely (q.s.) if there exists a polar set $D$ with ${\mathbb V}(D)=0$ such that it holds for each $\omega\in D^c$. We say that a property holds ${\cal P}$-q.s. means that it holds $P$-a.s. for each $P\in{\cal P}$. If an event $A$ fulfills $\mathbb V(A)=1$, then we call the event $A$ occurs $\mathbb V$-a.s.

\section{Delay-Dependent Asymptotic Stability of $G$-SDDEs}

  For the convenience of presentation, all processes take values in $\mathbb{R}$. Let  $\mathbb{R}_+ = [0,\infty)$.   If $A$ is a subset of $\Omega$, denote by $I_A$ its
indicator function.
Let $(\Omega , {\cal H}, \{{\Omega}_t\}_{t\ge 0}, \hat{\mathbb{E}}, \mathbb{V})$ be a generalized sublinear expectation space.  Let $(B(t))_{t\ge0}$
be one dimensional $G$-Brownian motion defined on the sublinear expectation space.

Let $f,g: \mathbb{R}^2 \times \mathbb{R}_+ \to \mathbb{R}$ be Borel measurable functions and $h: \mathbb{R}_+ \to \mathbb{R}$ be deterministic continuous function.  For $\tau>0,$ let $\de(t)\in[0,\tau],t\ge0$ satisfy $d\de(t)/dt=\dot{\de}(t)\le \bar\de$. Consider one dimensional highly nonlinear variable delay SDE driven by $G$-Brownian motion ($G$-SDDE)
\begin{align} \label{sdde1}
dX(t)= &f(X(t),X(t-\de(t)),t)dt\nonumber\\
&\quad + g(X(t),X(t-\de(t)),t)d<B(t)>+h(t)dB(t)
\end{align}
on $t\geq 0$ with nonrandom initial data
 \begin{equation} \label{2.2}
\{X(t):-\tau\leq t\leq 0\} = \eta\in C([-\tau,0]; \mathbb{R}).
\end{equation}

The uniqueness of solutions to  stochastic differential equations driven by $G$-Brownian motion has been proved under the coefficients satisfying non-Lipschitzian conditions, where the coefficients is often bounded by a linear function (see, e.g.,  Lin \cite{Lq2013}). In this paper, however, the coefficients of $G$-SDDE (\ref{sdde1}) cannot be bounded by a linear function. The coefficients in (\ref{sdde1}) are called highly nonlinear in terms of Fei {\sl et al.} \cite{FHMS} and Hu {\sl et al.} \cite{HMS}, the corresponding equations (\ref{sdde1}) are called highly nonlinear $G$-SDDEs. In Fei {\sl et al.} \cite{FFY}, the existence and uniqueness of solutions to $G$-SDDEs (\ref{sdde1}) are proved, while the stability and boundedness of solutions to $G$-SDDE (\ref{sdde1}) is investigated as well. We will consider highly nonlinear $G$-SDDEs which, in general, do not satisfy the linear growth condition in this paper.  Therefore, we impose the polynomial growth condition, instead of the linear growth condition.  Let us provide these conditions as an assumption for our aim.

\begin{assp} \label{A2.1}
 Assume that for any $b>0$, there exists a positive constant $K_b$ such that
 \begin{align*}
& |f(x,y,t)-f(\bar x,y,t)|\vee|g(x,y,t)-g(\bar x,y,t)|\\
&\quad \leq K_b(|x-\bar x|+|y-\bar y|)
 \end{align*}
for all $x, \bar x,y,\bar y\in \mathbb{R}^d$ with $|x|\vee |\bar x|\vee|y|\vee|\bar y|\leq b$ and all $t\in \mathbb{R}_+$.
 Assume moreover that there exist  three constants $K >0$, $q_1,q_2$ such that

\begin{align}\label{Pol}
 &|f(x,y,t)|\leq K(1+|x|^{q_1}+|y|^{q_1})\nonumber\\
  &|g(x,y,t)|\leq K(1+|x|^{q_2}+|y|^{q_2})\nonumber\\
  &|h(t)|\le K
\end{align}
 for all $x,y\in\mathbb{R}, t\in\mathbb{R}_+$.
\end{assp}

The condition (\ref{Pol}) with $q_1=q_2=1$ is the familiar linear growth condition. But, we
emphasise that we are here interested in highly nonlinear $G$-SDDEs which mean $q_1>1$ or $q_2>1$. The condition (\ref{Pol}) is referred as the polynomial growth condition.
   It is known that Assumption \ref{A2.1} only guarantees that the $G$-SDDE (\ref{sdde1}) with the initial data (\ref{2.2}) has a unique maximal solution,
which  may explode to infinity at a finite time.  To avoid such a possible explosion, we need to impose an
 additional condition by Lyapunov functions.  To this end, we need more notation.

 We denote $C^{2,1}(\mathbb{R} \times \mathbb{R}_+; \mathbb{R}_+)$ as the family of non-negative functions $U(x,t)$ defined on $(x,t)\in \mathbb{R} \times \mathbb{R}_+$ which are continuously twice differentiable in $x$ and once in $t$. Now we can state another assumption.

\begin{assp} \label{A2.2} Let $H(\cdot)\in C(\mathbb{R}\times[-\tau,\infty);\mathbb{R}_+)$.
Assume that there exists a function $\bar U\in C^{2,1}(\mathbb{R}\times \mathbb{R}_+; \mathbb{R}_+)$, nonnegative constants $c_1, c_2, c_3$ and $q=2(q_1\vee q_2$ ) such that
\begin{align}\label{H12}
c_3< c_2(1-\bar\de), |x|^q \leq \bar U(x,t) \leq H(x,t)
\end{align}
for $\forall(x,t) \in \mathbb{R}\times\mathbb{R}_+$, and
\begin{align}\label{ALV}
\mathbb{L}\bar U(x,y,t): =& \bar U_t(x,t)
  + \bar U_x(x,t) f(x,y,t)\nonumber\\
  &+  G(2\bar U_x(x,t)g(x,y,t)+\bar U_{xx}(x,t)h^2(t))\nonumber\\
     \le&  c_1 -c_2 H(x,t) + c_3 H(y,t-\de(t))
\end{align}
 for all $x,y\in\mathbb{R},t \in  \mathbb{R}_+$.
\end{assp}

 The following result gives the boundedness of the solution to $G$-SDDE (\ref{sdde1}) (see, e.g., Fei {\sl et al.} \cite[Theorem 5.2]{FFY}).

\begin{prop}\label{Pro}
Under Assumptions \ref{A2.1} and \ref{A2.2},  the highly nonlinear $G$-SDDE (\ref{sdde1}) with the initial data (\ref{2.2}) has unique global solution satisfying
\begin{align*}
\sup\limits_{-\tau\leq t<\infty}\hat{\mathbb{E}}|X(t)|^q<\infty.
\end{align*}
\end{prop}

Next, we will use the method of Lyapunov functionals to investigate the delay-dependent asymptotic stability. We define two  segments $\bar X(t): = \{X(t+s): -2\tau\leq s\leq 0\}$ for $t\geq 0$.
For $\bar X(t)$ to be well defined for $0\leq t < 2\tau$,
we set $X(s)=\eta(-\tau)$ for $s\in [-2\tau, -\tau)$.
We construct the Lyapunov functional as follows
\begin{align*}
&V(\bar X(t),t) = U(X(t),t)\\
 & +\theta \int_{-\tau}^0 \int_{t+s}^t \Big[ \tau |f(X(u),X(u-\de(u)),u)|^2\\
 &\quad+\tau\bar\sigma^2|g(X(u),X(u-\de(u),u))|^2+\bar\sigma^2|h(u)|^2 \Big] du ds
\end{align*}
for $t \geq 0$, where $U\in C^{2,1}(\mathbb{R}\times \mathbb{R}_+; \mathbb{R}_+)$
such that
\begin{align*}
\lim_{|x|\to\infty} \ [\inf_{t\in\mathbb{R}_+}\ U(x,t)]=\infty,
\end{align*}
  and $\theta$ is a positive number to be determined later while we set
\begin{align*}
&f(x,s)=f(x,0), \ g(x,s)=g(x,0)
\end{align*}
for all $x\in\mathbb{R}, s\in[-2\tau,\infty)$.
Applying the It\^o formula for $G$-Brownian motion (see \cite{P2010}) to $U(X(t),t)$, we get, for $t\ge0$,

 \begin{align*}
  &dU(X(t),t) \\
  &= \Big( U_t(X(t),t)+ U_x(X(t),t) f(X(t),X(t-\de(t)),t) \Big) dt   \\
& \quad+\Big(U_x(X(t),t)g(X(t),X(t-\de(t)),t)\\
&\quad+\frac12 h^2(t)U_{xx}(X(t),t)\Big)d<B(t)> \\
&\quad+ U_x(X(t),t)h(t)dB(t)\\
&\le {\LL} U(X(t),X(t-\de(t)),t)dt+U_x(X(t),t)h(t)dB(t)
\end{align*}
by \cite[Proposition 2.5]{FF}.
Rearranging terms gives
\begin{align*}
dU&(X(t),t)\\
 & \leq\Big( U_x(X(t),t)
 [ f(X(t),X(t-\de(t)),t) - f(X(t),X(t),t) ]  \\
 & +  { \mathcal{L}}U(X(t),X(t-\de(t)),t) \Big)dt+ U_x(X(t), t)h(t)dB(t),
\end{align*}
where the function ${\mathcal{ L}}U: \mathbb{R}^2 \times C([-\de(t),0];\mathbb{R}^d)\times \mathbb{R}_+ \to \mathbb{R}$  is defined by
\begin{align}\label{LU}
  {\mathcal{L}}U&(x,y,t)= U_t(x,t)
  +U_x(x,t) f(x,x,t)\nonumber\\
  &+  G(2g(x,y,t)U_{x}(x,t)+h(t)^2U_{xx}(x,t)).
\end{align}
 Moreover, the fundamental theory of calculus shows
\begin{align*}
& d\Big(\int_{-\tau}^0 \int_{t+s}^t  \Big[ \tau |f(X(u),X(u-\de(u)),u)|^2\\
&+\tau\bar\sigma^2|g(X(u),X(u-\de(u)),u)|^2+\bar\sigma^2|h(u)|^2 \Big] du ds\Big) \\
& =  \Big( \tau \Big[ \tau |f(X(t),X(t-\de(t)),t)|^2 \\
&+\tau\bar\sigma^2|g(X(t),X(t-\de(t)),t)|^2+\bar\sigma^2|h(t)|^2  \Big]\\
&-\int_{t-\tau}^t \Big[ \tau |f(X(u),X(u-\tau),u)|^2 \\
&+\tau\bar\sigma^2|g(X(u),X(u-\de(u)),u)|^2+\bar\sigma^2|h(u)|^2  \Big] du  \Big)   dt.
\end{align*}
\begin{lemma} \label{L3.1}
With the notation above, $V(\bar X(t), t)$ is $G$-It\^o process on $t\geq 0$ with its It\^o differential
\begin{equation*}
dV(\bar X(t), t) \leq LV(\bar X(t), t)dt + dM(t)
\quad\mbox{q.s.},
\end{equation*}
where  $M(t)$ is a $G$-continuous martingale with $M(0)=0$ and
\begin{align*}
  LV&(\bar X(t),t)  = \\
  & U_x(X(t),t) [ f((X(t),X(t-\de(t)),t)- f(X(t),X(t),t) ] \\
  & +  {\mathcal{L}}U(X(t),X(t-\de(t)),t)\\
& +\theta \tau \Big[ \tau |f(X(t),X(t-\de(t)),t)|^2 \\
&+\tau\bar\sigma^2|g(X(t),X(t-\de(t)),t)|^2+\bar\sigma^2|h(t)|^2  \Big]\\
&-\theta\int_{t-\tau}^t \Big[ \tau |f(X(u),X(u-\tau),u)|^2 \\
&+\tau\bar\sigma^2|g(X(u),X(u-\de(u)),u)|^2+\bar\sigma^2|h(u)|^2  \Big] du.
\end{align*}
\end{lemma}
To study the delay-dependent asymptotic stability of  the $G$-SDDE (\ref{sdde1}), we need to impose two new assumptions.
\begin{assp} \label{A3.2}
Assume that there are functions $U\in C^{2,1}(\mathbb{R} \times \mathbb{R}_+; \mathbb{R}_+)$, $U_1 \in C(\mathbb{R}\times[-\tau,\infty); \mathbb{R}_+)$,
and positive  numbers $\alpha_1,\alpha_2$ and $\beta_k$ ($k=1,2,3,4$) such that
\begin{equation}\label{A3.2c}
\alpha_2<\alpha_1(1-\bar\de)
\end{equation}
and
\begin{align}\label{A3.2c1}
&{\mathcal{L}}U(x,y,t) + \beta_1 |U_x(x,t)|^2\\
&+ \beta_2|f(x,y,t)|^2 +\beta_3|g(x,y,t)|^2+\beta_4|h(t)|^2 \nonumber\\
&\leq -\alpha_1 U_1(x,t)+\alpha_2U_1(y,t-\de(t)),
\end{align}
for all $x\in \mathbb{R}, t\in\mathbb{R}_+$.
\end{assp}
\begin{assp} \label{A3.3}
Assume that there exists a positive number $\varpi$ such that
\begin{equation*}
|f(x,x,t)-f(x,y,t)| \leq  \varpi|x-y|
\end{equation*}
for all $x\in\mathbb{R}, t\in [-2\tau,\infty)$.
\end{assp}

\begin{theorem}\label{T3.4}
Let Assumptions  \ref{A2.1}, \ref{A2.2}, \ref{A3.2} and \ref{A3.3} hold, where $h(\cdot)$ is a deterministic function.  Assume also that
\begin{equation}\label{WF}
 \tau\le \sqrt{\frac{4\beta_1\beta_2}{3\varpi^2}}\bigwedge\sqrt{\frac{4\beta_1\beta_3}{3\varpi^2\bar\sigma^2}}\bigwedge\frac{4\beta_1\beta_4}{3\varpi^2\bar\sigma^2}.
 \end{equation}
Then for any given initial data (\ref{2.2}), the solution of the $G$-SDDE (\ref{sdde1}) has the properties that
\begin{equation} \label{T3.5b}
\mathcal{E}\int_{0}^\infty  U_1(X(t),t)dt <\infty
\end{equation}
and
\begin{equation*}
\sup_{0\leq t <\infty} \hat{\mathbb{E}} U(X(t),t) < \infty.
\end{equation*}
\end{theorem}
\noindent
{\textbf{ Proof}}:   Fix the initial data $\eta\in C([-\tau,0];\mathbb{R})$ arbitrarily.
Let $k_0 >0$ be a sufficiently large integer such that $\|\eta\|:=\sup_{-\tau\leq s\leq0}|\eta(s)| < k_0$.
For each integer $k>k_0$, define the stopping time
$$
\sigma_k = \inf\{t \geq 0: |x(t)| \geq k\}.
$$
 It is easy to see that $\sigma_k$ is increasing as $k\to\infty$ and
$\lim_{k\to\infty} \sigma_k =\infty$ {\rm a.s.}    By the It\^o formula for $G$-Brownian motion
we obtain from Lemma \ref{L3.1} that
\begin{align} \label{intV}
 &\hat{\mathbb{E}} V(\bar X(t\wedge\sigma_k), t\wedge\sigma_k)\nonumber\\
&  \leq V(\bar X(0),0)+ \hat{\mathbb{E}} \int_0^{t\wedge\sigma_k} LV(\bar X(s), s) ds
\end{align}
for any $t\geq 0$ and $k\geq k_0$.
Let $\theta= 3\varpi^2/(4\beta_1) $.   By Assumption \ref{A3.3} and Cauchy-Schwartz inequality, it is easy to see that
\begin{align}
  &U_x(X(t),t)[f(X(t),X(t),t)- f(X(t),X(t-\de(t)),t) ]  \nonumber\\
  &\leq \beta_1|U_x(X(t),t) |^2 + \frac{\varpi^2}{4\beta_1}|X(t)-X(t-\de(t))|^2.
\end{align}
By condition (\ref{WF}), we also have
\begin{equation*}
  \theta\tau^2\le \beta_2,~\theta\bar\sigma^2\tau^2\le\beta_3,~\theta\bar\sigma^2\tau\leq \beta_4.
\end{equation*}
It then follows from Lemma \ref{L3.1}  that
\begin{align*}
  &L V(\bar X(s), s) \leq  {\mathcal{L}}U(X(s),X(s-\de(t)),s)
  + \beta_1 |U_x(X(s),s)|^2\\
  &  +\beta_2 |f(X(s),X(s-\de(s)),s)|^2\\
 & + \beta_3 |g(X(s),X(s-\de(s)),s)|^2+\beta_4|h(s)|^2\\
 & +  \frac{\varpi^2}{4\beta_1}|X(s)-X(s-\de(s))|^2 \\
 & - \frac{3\varpi^2}{4\beta_1}  \int_{s-\de(s)}^s
\Big[ \tau |f(X(u),X(u-\de(u)),u)|^2\\
&+\tau\bar\sigma^2|g(X(u),X(u-\de(u)),u)|^2+\bar\sigma^2|h(u)|^2 \Big] du.
\end{align*}
By Assumption \ref{A3.2c1}, we then have
\begin{align*}
  &L V (\bar X(s),s) \\
  & \leq -\alpha_1 U_1(X(s),s) +\alpha_2U_1(X(s-\de(s)),s-\de(s))\\
 & +  \frac{\varpi^2}{4\beta_1} \int_{-\de(s)}^0|X(s)-X(s+u)|^2du \\
 & - \frac{3\varpi^2}{4\beta_1}  \int_{s-\de(s)}^s
\Big[ \tau |f(X(u),X(u-\de(u)),u)|^2\\
&+\tau\bar\sigma^2|g(X(u),X(u-\de(u)),u)|^2+\bar\sigma^2|h(u)|^2 \Big] du.
\end{align*}
Substituting this into (\ref{intV}) implies
\begin{align} \label{intV2}
\mathbb{ E} V(\bar X(t\wedge\sigma_k),  t\wedge\sigma_k) \leq V(\bar X(0),0) +I_1+I_2,
\end{align}
where
\begin{align*}
&I_1 = \hat{\mathbb{E}}\int_0^{t\wedge\sigma_k} \big[ -\alpha_1 U_1(X(s),s)\\
&\quad +\alpha_2U_1(X(s-\de(s)),s-\de(s)) \big] ds, \\
&I_2 =  \frac{\varpi^2}{4\beta_1}  \hat{\mathbb{E}}\Big[\int_0^{t\wedge\sigma_k} \big[ |X(s)-X(s-\de(s))|^2 \\
&\quad- 3\int_{s-\de(s)}^s\big( \tau |f(X(u),X(u-\de(u)),u)|^2\\
&\quad+ \tau\bar\sigma^2|g(X(u),X(u-\de(u)),u)|^2+\bar\sigma^2|h(u)|^2\big) du\big ]  ds\Big].
\end{align*}
We notice the following fact
\begin{align*}
&\int_0^{\sigma_k\wedge t}U_1(X(s-\de(s)),s-\de(s))ds\nonumber\\
&\leq \frac{ 1}{1-\bar\de}\int_{-\tau}^0U_1(\eta(s),s)ds+\frac{1}{1-\bar\de}\int_0^{\sigma_k\wedge t}U_1(x(s),s)ds.
\end{align*}
Thus we get
\begin{align*}
&I_1\leq\frac{\alpha_2}{1-\bar\de}\int_{-\tau}^0  U_1(\eta(u),u)  du \\
&\quad+\frac{\bar\alpha}{1-\bar\de} \hat{\mathbb{E}}\Big[-\int_0^{t\wedge\sigma_k} U_1(X(s),s) ds\Big],
\end{align*}
where $\bar\alpha= (1-\bar\de)\alpha_1 - \alpha_2> 0$ by Assumption \ref{A3.2}.
Substituting this into (\ref{intV2}) yields
\begin{align} \label{intV3}
\frac{\bar\alpha}{1-\bar\de} \mathcal{E}\int_0^{t\wedge\sigma_k} U_1(X(s),s) ds \leq  C_1 +I_2,
\end{align}
where $C_1$ is a constant defined by
\begin{equation*}
C_1 = V(\bar X(0),0) + \frac{\alpha_2}{1-\bar\de} \int_{-\tau}^{0}  U_1(\eta(s),s)  ds.
\end{equation*}
Let $k \to\infty$ in (\ref{intV3}) to obtain
\begin{align} \label{intV4}
\frac{\bar\alpha}{1-\bar\de} \mathcal{E}\int_0^{t} U_1(X(s),s) ds \leq  C_1
 +\bar I_2,
\end{align}
where
\begin{align*}
\bar I_2 =&  \frac{\varpi^2}{4\beta_1}  \hat{\mathbb{E}}\Big[\int_0^t \big[ |X(s)-X(s-\de(s))|^2 \nonumber\\
&\quad- 3\int_{s-\de(s)}^s\big( \tau |f(X(u),X(u-\de(u)),u)|^2\\
&\quad+\tau\bar\sigma^2|g(X(u),X(u-\de(u)),u)|^2+\bar\sigma^2|h(u)|^2 \big) du\big ]  ds\Big].
\end{align*}
For $t\in [0,\tau]$, we  have
\begin{align}\label{H1H2}
\bar I_2   &\leq  \frac{\varpi^2}{2\beta_1}   \int_0^{\tau} (\hat{\mathbb{E}}|X(s)|^2 + \hat{\mathbb{E}}|X(s-\de(s))|^2 )ds\nonumber\\
&\leq \frac{\tau\varpi^2}{\beta_1}   \Big( \sup_{-\tau\leq u\leq\tau} \hat{\mathbb{E}}|X(u)|^2 \Big).
\end{align}
   For $t> \tau$, we have
\begin{align}\label{H2j}
 &\bar I_2  \le  \frac{\tau\varpi^2}{\beta_1} \Big( \sup_{-\tau\leq u\leq\tau} \hat{\mathbb{E}}|X(u)|^2 \Big)\nonumber\\
  &\quad+ \frac{\varpi^2}{4\beta_1}\hat{\mathbb{E}}\Big[\int_\tau^t \big[ |X(s)-X(s-\de(s))|^2 \nonumber\\
  &\quad- 3\int_{s-\de(s)}^s\big( \tau |f(X(u),X(u-\de(u)),u)|^2\nonumber\\
  &\quad+\tau\bar\sigma^2|g(X(u),X(u-\de(u)),u)|^2+\bar\sigma^2|h(u)|^2 \big) du\big ]  ds\Big].
\end{align}
 Noting
 \begin{align*}
&|X(s)-X(s-\de(s))|\\
 &    = \Big|\int_{s-\de(s)}^{s} f(X(u),X(u-\de(u)),u)du \\
 &+\int_{s-\de(s)}^{s} g(X(u),X(u-\de(u)),u)d<B(u)>\\
 &+\int_{s-\de(t)}^{s} h(u)dB(u)\Big|,
\end{align*}
by \cite[Proposition 2.5]{FF} and Cauchy-Schwartz inequality, we have
\begin{align*}
&|X(s) -X(s-\de(s))|^2 \\
 &\leq 3\int_{s-\de(s)}^{s}\tau |f(X(u),X(u-\de(t)),u)|^2du \\
 &+ 3\bar\sigma^2\int_{s-\de(s)}^{s}\tau |g(X(u),X(u-\de(s)),u)|^2du\\
 &+3\Big(\int_{s-\de(s)}^{s}h(u)dB(u)\Big)^2\quad \mbox{q.s.}
\end{align*}
Thus we get
\begin{align*}
&\hat{\mathbb{E}}\int_\tau^t \big[ |X(s)-X(s-\de(s))|^2 \\
&- 3\int_{s-\de(s)}^s\big( \tau |f(X(u),X(u-\de(u)),u)|^2\\
&+ \bar\sigma^2\tau|g(X(u),X(u-\de(u)),u)|^2+\bar\sigma^2|h(u)|^2\big) du\big ]  ds\\
 &\leq 3\hat{\mathbb{E}}\int_\tau^t\Big(\int_{s-\de(s)}^{s}h(u)dB(u)\Big)^2ds\\
 &- 3\bar\sigma^2\int_\tau^t  \int_{s-\de(s)}^s|h(u)|^2 du  ds=0
\end{align*}
by  Lemma \ref{BBI} and the properties of sublinear expectation and $G$-Brownian motion.

Thus from (\ref{H1H2}) and (\ref{H2j}) we get
\begin{align}\label{H3H4}
\bar I_2 \leq\frac{\tau\varpi^2}{\beta_1}   \Big( \sup_{-\tau\leq u\leq\tau} \hat{\mathbb{E}}|X(u)|^2 \Big).
\end{align}
Substituting (\ref{H3H4}) into (\ref{intV4}), together with (\ref{WF}), yields
$$\frac{\bar\alpha}{1-\bar\de}\mathcal{E}\int_0^{t} U_1(X(s),s) ds  \leq C_1+ \frac{4\beta_3}{3\bar\sigma^2} \sup_{-\tau\leq u\leq\tau} \hat{\mathbb{E}}|X(u)|^2 :=C_2. $$
Letting $t\to \infty$ gives
\begin{align*}
\mathcal{E}\int_{0}^{\infty}  U_1(X(s),s) ds
 \leq \frac {(1-\bar\de) C_2}{\bar\alpha}.
\end{align*}
We deduce easily from (\ref{intV2}) that
\begin{align*}
\hat{\mathbb{E}}U \Big( X(t\wedge\sigma_k),
t\wedge\sigma_k \Big)
 \leq C_1  +I_2.
\end{align*}
Letting $k\to \infty$ we get
\begin{align*}
\hat{\mathbb{E}}U (X(t),t) \leq C_2 <\infty  ,
\end{align*}
which shows
\begin{align*}
\sup_{0\leq t <\infty}\hat{\mathbb{E}}U(X(t),t)    <\infty.
\end{align*}
Thus the proof is complete. $\Box$
\begin{coro} \label{C3.5}
Let the conditions of Theorem \ref{T3.4} hold.   If there moreover exists a  pair of positive constants $c$  and $p$ such that
\begin{align*}
  c  |x|^p \leq U_1(x), \quad \forall x\in \mathbb{R} \times \mathbb{R}_+,
\end{align*}
then for any given initial data (\ref{2.2}), the solution of $G$-SDDE (\ref{sdde1})  satisfies
\begin{align*}
 {\mathcal E}\Big [ \int_0^{\infty} |X(t)|^p dt\Big]  < \infty,
\end{align*}
which means
\begin{align}\label{C3.5b}
\int_0^{\infty} {\mathcal E}|X(t)|^p dt  < \infty.
\end{align}
\end{coro}

Noting the property of sublinear $\cal E$, this corollary follows from Theorem \ref{T3.4} obviously. However, it does not follow from (\ref{C3.5b}) that
$\lim_{t\to\infty}{\mathcal E}|X(t)|^p=0$. The following proposition provides a stronger claim.
\begin{prop}\label{T3.6}
Let the conditions of Corollary  \ref{C3.5} hold. If, moreover,
$$
p\geq 2 \quad\hbox{and}\quad  (p+q_1-1)\vee(p+q_2-1)\leq  q,
$$
then the solution of $G$-SDDE (\ref{sdde1}) satisfies
$$
 \lim_{t\to\infty}\hat{\mathbb{E}}|X(t)|^p=0
$$
for any initial data (\ref{2.2}).  That is, $G$-SDDE (\ref{sdde1}) is  p-moment asymptotically stable under the sublinear expectation $\hat{\mathbb E}$.
\end{prop}

\noindent
\textbf{Proof}:   Fix the initial data (\ref{2.2}) arbitrarily.  For any $0\leq t_1 < t_2<\infty$,
by the It\^o formula for $G$-Brownian motion, we get
\begin{align*}
&d|X(t)|^p=p|X(t)|^{p-1}f(X(t),X(t-\de(t)),t)dt \\
&+ p|X(t)|^{p-2}[X(t)g(X(t),X(t-\de(t)))\\
&+\frac{1}{2}(p-1)h^2(t)]d<B(t)>+p|X(t)|^{p-1}h(t)dB(t),
\end{align*}

Noting $d<B(t)>\le \bar\sigma^2dt~$q.s. (see, e.g., Fei and Fei \cite{FF1}), along with the polynomial growth condition (\ref{Pol}), we get
\begin{align*}
&\big | \hat{\mathbb{E}}|X(t_2)|^p-\hat{\mathbb{E}}|X(t_1)|^p \big|\\
&\leq \hat{\mathbb{E}} \int_{t_1}^{t_2} \Big( p|X(t)|^{p-1} |f(X(t),X(t-\de(t)),t)|\\
& \quad+ p\bar\sigma^2 |X(t)|^{p-1} |g(X(t),X(t-\de(t)),t)|\\
&\quad+\frac{1}{2}p(p-1)\bar\sigma^2|X(t)|^{p-2}h^2(t)\\
& \leq \hat{\mathbb{E}} \int_{t_1}^{t_2} \Big( pK|X(t)|^{p-1}  \big[1+|X(t)|^{q_1} +|X(t-\de(t))|^{q_1}\big]\\
&\quad +pK\bar\sigma^2|X(t)|^{p-1}  \big[1+|X(t)|^{q_2} +|X(t-\de(t))|^{q_2}\big]\\
&\quad+ \frac{1}{2}p(p-1)K^2\bar\sigma^2 |X(t)|^{p-2}  \Big)dt.
\end{align*}
By inequalities,
\begin{align*}
&|X(t)|^{p-1}  |X(t-\de(t))|^{q_1} \\
&\quad\leq |X(t)|^{p+q_1-1} + |X(t-\de(t))|^{p+q_1-1} ,\\
&|X(t)|^{p-1} \leq 1 + |X(t)|^q,
\end{align*}
etc., and noting that for any $1\leq \bar p\leq q$, by Proposition \ref{Pro} we have
\begin{align*}
\hat{\mathbb{E}}|X(t-\de(t))|^{\bar p}\leq 1+\sup\limits_{-\tau\leq s< \infty}\hat{\mathbb{E}}|X(s)|^q<\infty,
\end{align*}
we can obtain
 \begin{align*}
 \big| \hat{\mathbb{E}}|X(t_2)|^p-&\hat{\mathbb{E}}|X(t_1)|^p \big|
 \leq C_3(t_2-t_1),
\end{align*}
where
\begin{align*}
&C_3=(4pK(1+\bar\sigma^2)\\
&\quad+\frac{1}{2}K^2p(p-1)\bar\sigma^2)(1+\sup_{-\tau \leq s< \infty}\hat{\mathbb{E}}|X(s)|^q)<\infty.
\end{align*}
Thus we have $\hat{\mathbb{E}}|X(t)|^p $ is uniformly continuous in $t$ on $\mathbb{R}_+$.
By (\ref{C3.5b}), there is a sequence $\{t_l\}_{l=1}^\infty$ in $\mathbb{R}$ such that
${\mathbb E}|X(t_l)|^p\rightarrow0,$ which easily show  the claim. Thus, the proof is complete. $\Box$

\begin{prop}\label{P3.7}
Let the conditions of Theorem \ref{T3.4} hold.  Assume also that there are positive constants $p$ and $c$ such that
 \begin{align} \label{P3.7a}
 c|x|^p\leq U(x,t), \quad \forall (x,t)\in \mathbb{R}\times\mathbb{R}_+ .
 \end{align}
 Moreover assume there exists a function $W: \mathbb{R}\rightarrow \mathbb{R}_+$ such that
  \begin{align*}
  W(x)=0 ~\mbox{if and only if}~ x=0
   \end{align*}
 and
 $$W(x)\leq U_1(x),\forall x\in\mathbb{R}.$$
  Then for any given initial data (\ref{2.2}), the solution $X(\cdot)$ to Eq. (\ref{sdde1}) obeys that
 \begin{align*}
\lim_{t\to\infty}X(t)=0 \quad \mathbb{V}\mbox{-a.s.}
\end{align*}
\end{prop}

\noindent\textbf{Proof}:
  Let $X(\cdot)$ be the solution to Eq. (\ref{sdde1}) with initial data $\eta$ defined  in (\ref{2.2}). Since the conditions in Theorem \ref{T3.4} hold, we can show that
$$
C_4:={\mathcal E}\Big[\int_{0}^\infty  W(X(t))dt\Big] <\infty,
$$
which implies
\begin{align}\label{Wx}
\int_{0}^\infty  W(X(t))dt < \infty \quad \mathbb{V}\mbox{-a.s.}
\end{align}
Set $\sigma_k:=\inf\{t\geq0: |X(t)|\ge k\}$.
We observe from (\ref{Wx}) that
\begin{align}\label{Wxt}
\lim_{t\to\infty}\inf W(X(t))=0 \quad\mathbb{V}\mbox{-a.s.}
\end{align}
 Moreover,  in the same way as Theorem \ref{T3.4} was proved, we can show that
\begin{align*}
\hat{\mathbb{E}}|X(T\wedge\sigma_k)|^p &\leq C, \quad\forall T >0,
\end{align*}
which implies, by Chebyshev inequality (see, e.g., Chen {\sl et al.} \cite{CWL}, Proposition 2.1 (2)),
$$
k^p \mathbb{V}(\sigma_k\leq T) \leq C.
$$
Letting $T\to\infty$ yields
 \begin{align} \label{P5.6b}
 k^p\mathbb{V}(\sigma_k< \infty) \leq C.
 \end{align}

We now claim that
\begin{align}\label{P5.6.1}
\lim_{t\to\infty} W(X(t))=0 \quad \mathbb{V}\mbox{-a.s.}
\end{align}
In fact, if this is false, then we can find a  number $\varepsilon \in (0,1/4)$ such that
\begin{align}\label{P5.6.2}
\mathbb{V} (\Omega_1)\geq 4\varepsilon,
\end{align}
where $\Omega_1=\{\lim\sup _{t\to\infty} W(X(t))>2\varepsilon\}$.
 Recalling (\ref{P5.6b}), we can find an integer $m$ sufficiently large for
${\cal V}(\sigma_{m} < \infty)\leq\mathbb{V}(\sigma_{m} < \infty) \leq \varepsilon$.  This means that
\begin{align}\label{P5.6.3}
\mathbb{V}(\Omega_2)=1-\mathcal{V}(\Omega_2^c) \geq 1-\varepsilon.
\end{align}
where $\Omega_2:=\{|X(t)|< m \hbox{ for }\forall t\geq -\tau\}$, and  $\Omega_2^c$ is the complement of $\Omega_2$.
By (\ref{P5.6.2}) and (\ref{P5.6.3}) we get
\begin{align}\label{P5.6.4}
\mathbb{V}(\Omega_1\cap\Omega_2) \ge  \mathbb{V}(\Omega_1) - \mathcal{V}(\Omega_2^c) \geq 3\varepsilon.
\end{align}

Let us now define the stopped process $\zeta(t)=X(t\wedge\sigma_m)$ for $t\geq -\tau$.
Clearly, $\zeta(t)$ is a bounded It\^o process with its differential
\begin{align*}
d\zeta(t) = \phi(t)dt +\psi(t)d<B(t)>+h(t)dB(t),
\end{align*}
where
\begin{align*}
&\phi(t)=f(X(t),X(t-\de(t)),t)I_{[0,\sigma_m)}(t),\\
&\psi(t)=g(X(t),X(t-\de(t)),t)I_{[0,\sigma_m)}(t),
\end{align*}
where $f,g,h$ are defined by (\ref{sdde1}). Recalling the polynomial growth condition (\ref{Pol}) and the boundedness of $h(\cdot)$, we know that
$\phi(\cdot),\psi(\cdot)$ and $h(\cdot)$ are bounded processes,  say
\begin{align}\label{P5.6.5}
|\phi(t)|\vee |\psi(t)|\vee h(t)\leq C_{5} \quad\mathbb{V}\mbox{-a.s.}
\end{align}
for all $t\geq 0$ and some $C_5>0$.    Moreover, we also observe that $|\zeta(t)|\leq m$ for all $t\geq -\tau$.
Define a sequence of stopping times
$$ \rho_1= \inf\{t\geq 0:W(\zeta(t)) \geq 2\varepsilon\},$$
$$ \rho_{2j} = \inf\{t\geq \rho_{2j-1}:W(\zeta(t)) \leq \varepsilon\}, \quad j=1,2,\cdots,$$
$$ \rho_{2j+1}= \inf\{t\geq \rho_{2j}:W(\zeta(t)) \geq 2\varepsilon\}, \quad j=1,2,\cdots.$$
From (\ref{Wxt}) and the definition of $\Omega_1$ and $\Omega_2$, we have
\begin{align}\label{OO}
\Omega_1\cap\Omega_2 \subset \{\sigma_m = \infty\} \bigcap \Big(  \cap_{j=1}^\infty \{\rho_j < \infty\} \Big).
\end{align}
We also note that for all $\omega\in\Omega_1\cap\Omega_2$,  and $j\geq 1$,
 \begin{align}\label{P5.6.6}
W(\zeta(\rho_{2j-1}))-W(\zeta(\rho_{2j})) =\varepsilon \quad \hbox{ and }  \nonumber \\
W(\zeta(t))\geq \varepsilon \quad \hbox{ when }    t\in [\rho_{2j-1}, \rho_{2j}].
\end{align}
Since $W(\cdot)$ is  uniformly continuous in the close ball $\bar S_m=\{x \in \mathbb{R}^d:|x|\leq m \}$. We can choose $\delta=\delta(\varepsilon)>0$ small sufficiently for which
\begin{align} \label{Wuv}
 |W(\zeta)-W(\bar\zeta)|<\varepsilon,
  \zeta, \bar\zeta \in \bar S_m,   \hbox{ with}  \quad  |\zeta-\bar\zeta|<\delta .
\end{align}
We emphasize that for $\omega\in\Omega_1\cap\Omega_2$, if $|\zeta(\rho_{2j-1}+u)-\zeta(\rho_{2j-1})| <\delta$ for all $u\in [0,\lambda]$ and some $\lambda >0$, then $\rho_{2j}-\rho_{2j-1} \geq \lambda$.
 Choose a sufficiently small positive number $\lambda$ and then a sufficiently large positive integer $j_0$ such that
\begin{align}\label{lambdadef}
 3 C_5^2  \lambda (\lambda+ \lambda C_1(2,{\bar\sigma})+C_2(2,{\bar\sigma}))\leq  \varepsilon \delta^2  \ \hbox{ and }  \  C_4< \varepsilon^2 \lambda  j_0.
\end{align}
By (\ref{P5.6.4}) and (\ref{OO}) we can further choose a sufficiently large number $T$ for
\begin{align}\label{P5.6.7}
\mathbb{V}(\rho_{2j_0}\leq T) \geq 2\varepsilon.
\end{align}
In particular, if $\rho_{2j_0}\leq T$, then $|\zeta(\rho_{2j_0})|<m$, and hence $\rho_{2j_0} < \sigma_m$ by the definition of $\zeta(t)$.   We hence have
\begin{align*}
\zeta(t,\omega)=X(t,\omega) \hbox{ for all $0\leq t\leq \rho_{2j_0}$ and $\omega\in \{\rho_{2j_0}\leq T\}$. }
\end{align*}
By the Burkholder-Davis-Gundy inequality under sublinear expectation (see, e.g., Lemma \ref{BDG}),  we can have that, for $1\leq j\leq j_0$,
\begin{align*}
  \hat{\mathbb{E}}\Big( & \sup_{0\leq t\leq \lambda}  |\zeta(\rho_{2j-1}\wedge T+t) -\zeta(\rho_{2j-1}\wedge T)|^2 \Big) \\
\leq & 3\lambda \hat{\mathbb{E}}\int_{\rho_{2j-1}\wedge T}^{\rho_{2j-1}\wedge T+\lambda} |\phi(s)|^2 ds \\
&\quad+3C_1(2,{\bar\sigma})\lambda \hat{\mathbb{E}}\int_{\rho_{2j-1}\wedge T}^{\rho_{2j-1}\wedge T+\lambda} |\psi(s)|^2 ds\\
 &\quad+ 3C_2(2,{\bar\sigma}) \int_{\rho_{2j-1}\wedge T}^{\rho_{2j-1}\wedge T+\lambda} |h(s)|^2ds  \\
 \leq &  3 C_5^2  \lambda (\lambda+ \lambda C_1(2,{\bar\sigma})+C_2(2,{\bar\sigma})),
\end{align*}
which, together with (\ref{lambdadef}) and Chebyshev inequality  for sublinear expectation $\hat{\mathbb{E}}$,  we can obtain that
$$
\mathbb{V}\Big( \sup_{0\leq t\leq \lambda}  | \zeta(\rho_{2j-1}\wedge T+t) -\zeta(\rho_{2j-1}\wedge T)| \geq \delta \Big) \leq \varepsilon.
 $$
 Noting that $\rho_{2j-1}\le T$ if $\rho_{2j_0}\le T$, we
can derive from (\ref{P5.6.7}) and the above inequality  that
\begin{align*}
 \mathbb{V}&\Big(\{\rho_{2j_0}\leq T\} \cap \Big\{\sup_{0\leq t\leq \lambda}
 |\zeta(\rho_{2j-1}+t)  - \zeta(\rho_{2j-1})|  < \delta \Big\} \Big) \\
  &=   \mathbb{V}(\rho_{2j_0}\leq T) \\
  &   -\mathcal{V}\Big( \{\rho_{2j_0}\leq T\}  \cap \Big\{\sup_{0\leq t\leq \lambda}
 |\zeta(\rho_{2j-1}+t)  - \zeta(\rho_{2j-1})| \geq \delta \Big\} \Big) \\
   &\geq   \mathbb{V}(\rho_{2j_0}\leq T) -\mathcal{V}\Big( \sup_{0\leq t\leq \lambda}   |\zeta(\rho_{2j-1}+t)  - \zeta(\rho_{2j-1})| \geq \delta \Big) \\
  &\geq  \varepsilon.
\end{align*}
This, together with (\ref{Wuv}),   implies easily that
\begin{align} \label{P5.6.12}
\mathbb{V}\Big( \{\rho_{2j_0}\leq T\} \cap \{\rho_{2j}-\rho_{2j-1} \geq \lambda\} \Big) \geq \varepsilon.
\end{align}
 By (\ref{P5.6.6}), (\ref{P5.6.12}), and the sub-additivity and Chebyshev inequality for sublinear expectation $\mathcal E$, we derive
\begin{align*}
 C_4 &   \geq \sum_{j=1}^{j_0}  {\mathcal E} \Big(  I_{\{\rho_{2j_0}\leq T\}} \int_{\rho_{2j-1}}^{\rho_{2j} }  W(X(t)) dt \Big)
\\
& \geq \varepsilon \sum_{j=1}^{j_0} {\mathcal E}  \Big( I_{\{\rho_{2j_0}\leq T\}} (\rho_{2j}-\rho_{2j-1})  \Big) \\
& \geq \varepsilon \lambda \sum_{j=1}^{j_0} \mathcal{V}\Big( \{\rho_{2j_0}\leq T\} \cap \{\rho_{2j}-\rho_{2j-1} \geq \lambda\}   \Big) \\
& \geq \varepsilon^2 \lambda j_0.
\end{align*}
This contradicts the second inequality in (\ref{lambdadef}).  Thus (\ref{P5.6.1}) must hold.

We now claim $\lim_{t\to \infty} X(t)=0$ $\mathbb{V}\mbox{-a.s.} $.  If this were not true, then
$$
\varepsilon_1 :={\cal V}(\Omega_3) > 0,
$$
where $\Omega_3 = \{\limsup_{t\to\infty} |X(t)| > 0\}$.   On the other hand, by (\ref{P5.6b}), we can find a positive integer $m_0$ large enough for
$\mathbb{V}(\sigma_{m_0} <\infty) \le 0.5 \varepsilon_1$.   Let $\Omega_4 = \{\sigma_{m_0}=\infty\}$.  Then
$$
{\cal V}(\Omega_3\cap \Omega_4)={\cal V}(\Omega_3)-\mathbb{V}(\Omega_3\cap\Omega_4^c) \geq {\cal V}(\Omega_3) - \mathbb{V}(\Omega_4^c) \ge 0.5\varepsilon_1.
$$
For any $\omega\in \Omega_3\cap \Omega_4$, $X(t,\omega)$ is bounded on $t\in \mathbb{R}_+$.  We can then find a sequence $\{t_j\}_{j\geq 1}$
such that  $t_j\to\infty$ and $X(t_j,\omega) \to \tilde X(\omega) \not= 0$ as $j\to \infty$. This, together with the continuity of $W$,  implies
$$
\lim_{j\to\infty} W(X(t_j,\omega)) = W(\tilde X(\omega)) > 0,
$$
which show
$$
\limsup_{t\to\infty} W(X(t,\omega)) > 0 \  \hbox{ for all $\omega\in \Omega_3\cap \Omega_4$.}
$$
But this contradicts (\ref{P5.6.1}).  We therefore must have the assertion $\lim_{t\to \infty} X(t)=0\quad\mathbb{V}\mbox{-a.s.}$ Thus, the proof is complete. $\Box$

\section{ Stability Example for $G$-SDDEs}

For a highly nonlinear $G$-SDDE (\ref{sdde1}) with $g(t)\equiv0$, by the equal-spaced partition of time, it follows that
\begin{align}\label{Dis}
X(t_{i})=&X(t_{i-1})+f(X(t_{i-1}),X(t_{i-1}-\de(t_{i-1})),t_{i-1})\Delta t_i\nonumber\\
&\quad+h(t_{i-1})(B(t_{i})-B(t_{i-1}))
\end{align}
with initial data $\eta=\{X(t)=\eta(t), t\in[-\tau,0]\},\Delta =\Delta {t_i}=t_i-t_{i-1}=1/N, i=1,\cdots,N$, where $N$ is a positive integer.

Let us first introduce the simulation algorithm for $G$-Brownian motion $(B(t))_{t\ge0}$. Now consider a random variable $\zeta=B({t_i})-B({t_{i-1}})\sim {\cal N}(0,[\underline\sigma^2,\bar\sigma^2]\Delta), i=1,\cdots,N$, we construct an experiment as follows. We take equal-step points $\sigma^k,k=1,\cdots,m,$ such that $\underline\sigma=\sigma^1<\sigma^2\cdots<\sigma^k<\cdots<\sigma^m=\bar\sigma$. For the $k$th-random sampling ($k=1,\cdots,m$),  $\zeta^{kj}_i (i=1,\cdots,N; j=1,\cdots,n)$ are from the classical normal distribution ${\cal N}(0, (\sigma^k)^2\Delta)$. Define $X^{kj}(t_i)$ by
\begin{align}\label{Dis1}
&X^{kj}(t_{i})=X^{kj}(t_{i-1})\nonumber\\
&+f(X^{kj}(t_{i-1}),X^{kj}(t_{i-1}-\de(t_{i-1})),t_{i-1})\Delta +h(t_{i-1})\zeta^{kj}_i
\end{align}
for $i=1,\cdots,N; j=1,\cdots,n; k=1,\cdots,m$.
 By the $\varphi$-max-mean algorithm from Peng \cite[Section 3.3]{P-scm}, for $p\geq 1$ we can obtain the estimator $\hat{\bar\mu}_{m,n}(p)$ of $\hat{\mathbb E}|X({t_i})|^p$ equals to $\max\limits_{1\leq j\leq n}\{\frac{1}{m}\sum\limits_{k=1}^m|X^{kj}({t_i})|^p\}$, and the estimator $\hat{\underline\mu}_{m,n}(p)$ of ${\mathcal E}|X({t_i})|^p$ to $\min\limits_{1\leq j\leq n}\{\frac{1}{m}\sum\limits_{k=1}^m|X^{kj}({t_i})|^p\}, i=1,\cdots,N$, respectively.

Also, in terms of Theorem 24 in Peng \cite{P-scm}, we know $\hat{\bar\mu}_{m,n}(p)=\max\limits_{1\leq j\leq n}\frac{1}{m}\sum\limits_{k=1}^m|X^{kj}({t_i})|^p$ is the unbiased estimator of $\hat{\mathbb E}|X({t_i})|^p$, and $\hat{\underline\mu}_{m,n}(p)=\min\limits_{1\leq j\leq n}\frac{1}{m}\sum\limits_{k=1}^m|X^{kj}({t_i})|^p$ is the unbiased estimator of $\mathcal E|X({t_i})|^p$. Moreover, we easily know that if $\hat{\mathbb E}|X(t)|$ with the solution $X(t)$ to highly nonlinear $G$-SDDE (\ref{sdde1}) is stable, then the solution $X(t)$ to $G$-SDDE (\ref{sdde1}) is qusi surely stable, i.e., ${\mathcal V}(\lim\limits_{t\rightarrow\infty}|X(t)|=0)=1$. Meanwhile, if ${\mathcal E}|X(t)|$ is unstable under sublinear expectation, then the solution $X(t)$ to $G$-SDDE (\ref{sdde1}) is unstable, i.e., ${\mathbb V}(\lim\sup\limits_{t\rightarrow\infty}|X(t)|\neq0)=1$.

Through utilizing the above discussion on the simulation algorithm, we now investigate an example to illustrate our theory.
\begin{expl} \label{E4.1}
{\rm  Let us consider the $G$-SDDE (\ref{$G$-SDDE1}), where $f(x,y,t)=-x^3-y,h(t)=\frac{1}{2}e^{-t}, \underline{\sigma}^2=0.5,\bar\sigma^2=1$ and $\bar\de=0.1$. We consider two case: $\de(t)=0.01$ and $\de(t)=2$ for all $t\geq 0$. In case $\de(t)=0.01$, let the initial data $\eta(u)=2+\text{sin}(u)$ for $u\in [-0.01,0],$  the sample paths of the solution of the $G$-SDDE are shown in Figure 4.1, which indicates that the $G$-SDDE is asymptotically stable. In the case $\de(t)=2$, let the initial data $\eta(u)=2+\text{sin}(u)$ for $u\in [-2,0],$ the solution of the $G$-SDDE (\ref{$G$-SDDE1}) are plotted in Figure 4.2, which indicates that $G$-SDDE is asymptotically unstable. From the example we can see $G$-SDDE (\ref{$G$-SDDE1}) is stable or not depends on how long or short the time-delay is.
 \begin{figure}[!hbt]
\begin{center}
\includegraphics[angle=0,height=8cm,width=9cm]{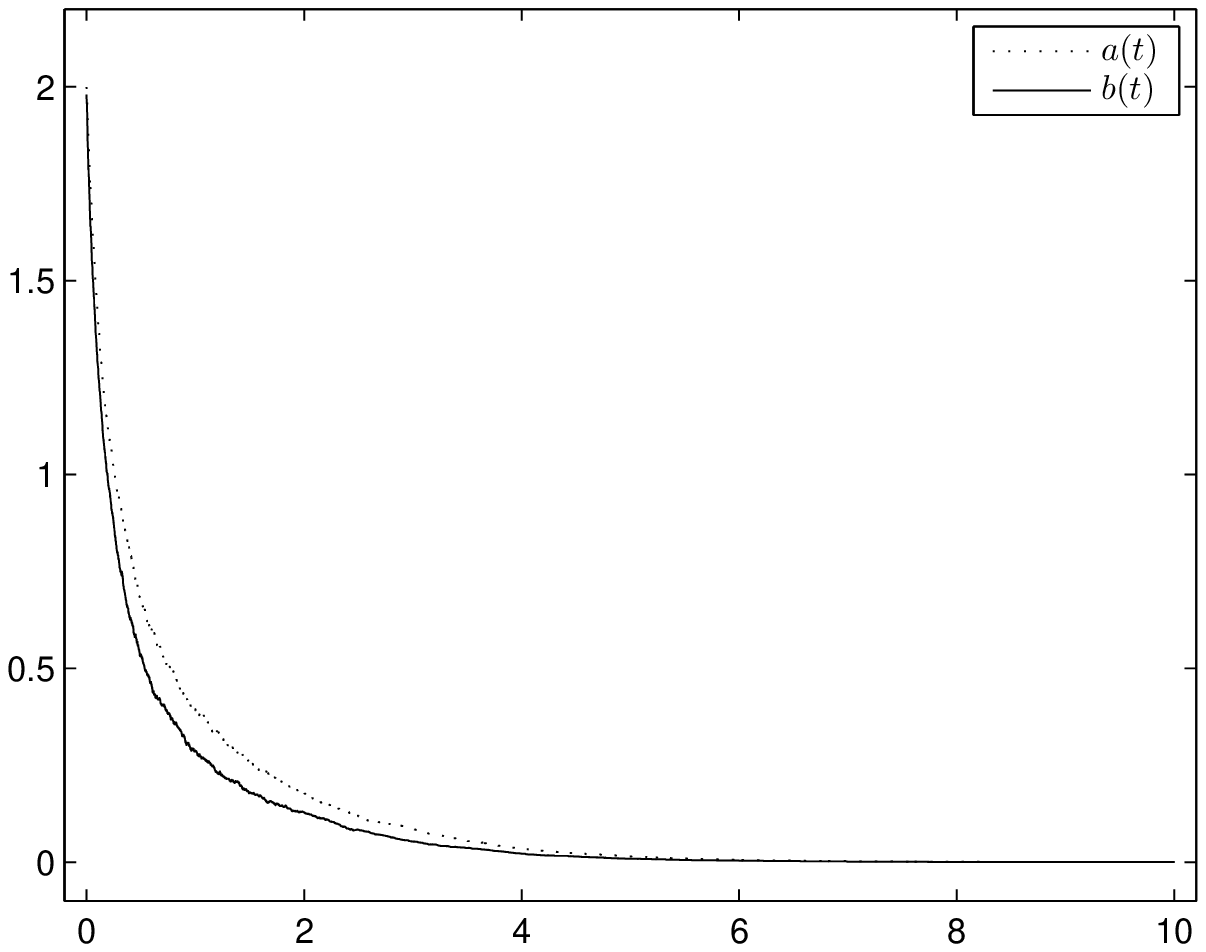}
\end{center}
\begin{center}
\small{Figure 4.1: The computer simulation of the upper expectation
$a(t)=\hat{\mathbb E}|X(t)|$ and the lower expectation $b(t)={\cal  E}|X(t)|$  of the solution to the $G$-SDDE  (\ref{$G$-SDDE1}) with $\de(t)=0.01$
using the Euler--Maruyama method with step size $10^{-3}$.}
\end{center}
\end{figure}

\begin{figure}[!h]
\begin{center}
\includegraphics[angle=0,height=8cm,width=9cm]{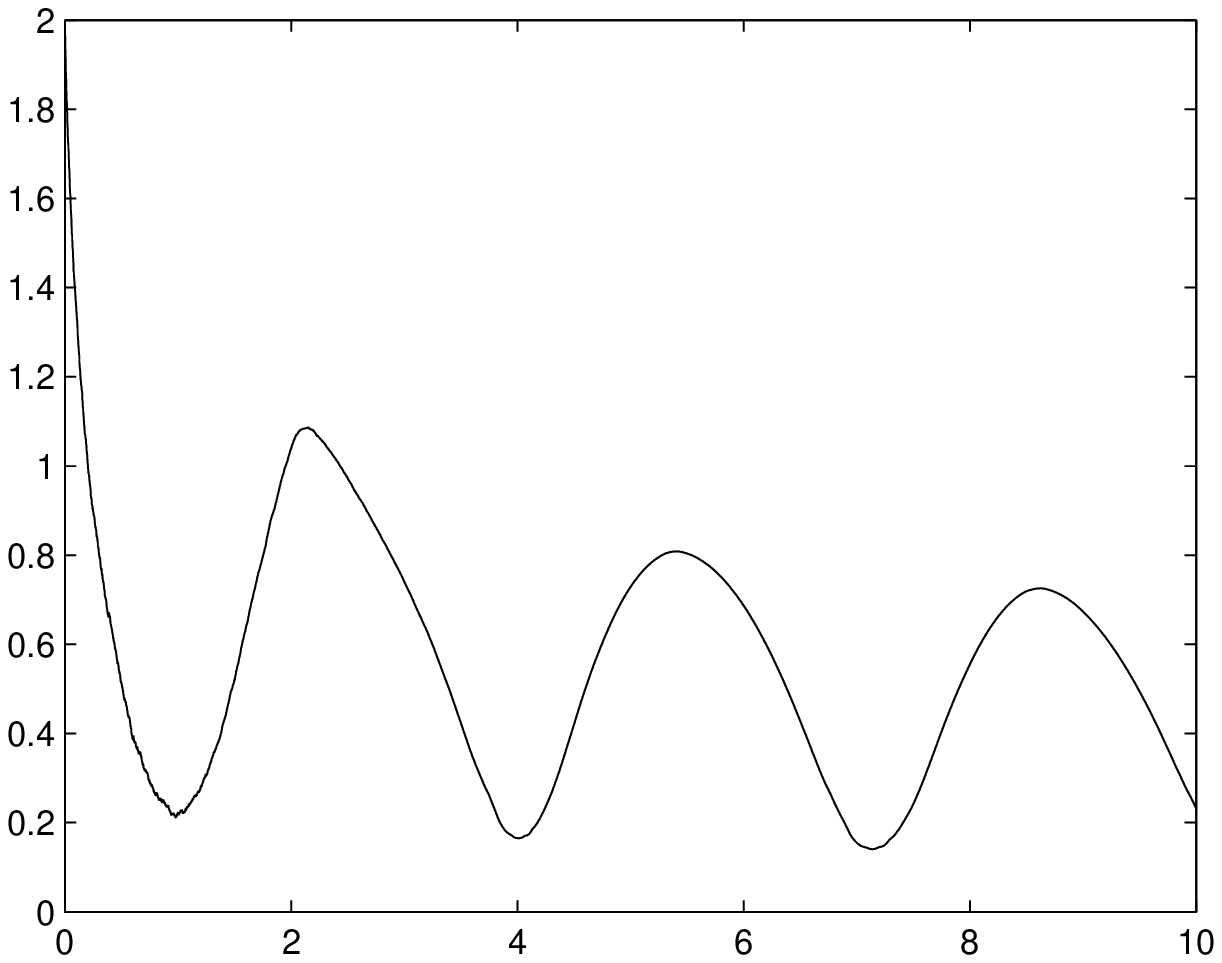}
\end{center}
\begin{center}
\small {Figure 4.2
: The computer simulation of
the lower expectation ${\cal  E}|X(t)|$  of the solution to the $G$-SDDE  (\ref{$G$-SDDE1}) with $\de(t)=2$
using the Euler--Maruyama method with step size $10^{-3}$.}
\end{center}
\end{figure}

 We  can see coefficients defined by (\ref{$G$-SDDE1})  satisfy Assumption \ref{A2.1} with $q_1=3$ and $q_2=0$. Define $\bar U(x,t)=|x|^6$  for  $(x,t)\in \mathbb{R}\times\mathbb{R}_+$. It is easy to show that
\begin{align*}
\mathbb{L}\bar U(x,y,t)= 6x^5f(x,y,t)+15x^4|h(t)|^2
\end{align*}
for $(x,t)\in \mathbb{R}\times\mathbb{R}_+$. Thus we get
\begin{align*}
&\mathbb{L}\bar U(x,y,t)= 6x^5(-x^3-y)+\frac {15}{4}x^4e^{-t}\\
&\leq 4x^4+5x^6-6x^8+y^6\\
&\leq c_1-2(1+1.5x^6+2.5x^8)+(1+1.5y^6+2.5y^8),
\end{align*}
where
\begin{align*}
c_1=\sup\limits_{x\in\mathbb{R}}(1+4x^4+8x^6-x^8)>0.
\end{align*}
 Thus, we can set $H(x,t)=1+1.5x^6+x^8$. Due to $1=c_3<(1-\bar\de)c_2=1.8$, we know Assumption \ref{A2.2} holds. From Proposition \ref{Pro}, the solution of the $G$-SDDE (\ref{$G$-SDDE1}) satisfies
 $$ \sup_{-\tau\leq t <\infty}\hat{\mathbb{E}}|X(t)|^6<\infty.$$
 To verify Assumption \ref{A3.2}, we define
 \begin{align}\label{UF}
U(x,t) =e^{-t}+x^2+ x^4,
\end{align}
which shows
\begin{align*}
U_t=-e^{-t},U_x(x,t) =2x+4x^3,U_{xx}=2+12x^2,
\end{align*}
for $(x,t)\in \mathbb{R}\times\mathbb{R}_+$. By equation (\ref{LU}), we have
\begin{align*}
&{\cal L}U(x,y,t)=\\
&-e^{-t}+(2x+4x^3)(-x^3-x)+\frac{1}{8}e^{-2t}(2+12x^2)\\
&\leq-\frac{3}{4}e^{-t}-0.5x^2-6x^4-4x^6.
\end{align*}
 Moreover
\begin{align}\label{Ui}
|U_{x}(x,t)|^2 = 4x^2+16x^4+16x^6,
\end{align}
\begin{align}\label {Fi2}
|f(x,y,t)|^2 =|-x^3-y|^2 \leq 1.5x^4+x^6+y^2+0.5y^4.
\end{align}
  Since $g(\cdot)\equiv0$, $\beta_3$ can take an arbitrary large number. Next, setting $\beta_1=0.1$, $\beta_2=0.05$, $\beta_4=1$, we get that
 \begin{align*}
 &{\cal L}U(x,y,t) +\beta_1 |U_x(x,t)|^2 + \beta_2|f(x,y,t)|^2 + \beta_4|h(t)|^2\\
&\leq -0.5e^{-t}-0.1x^2-4x^4-2x^6+0.05y^2+0.0025y^4\\
&\leq -(0.5e^{-t}+0.1x^2+4x^4+2x^6)\\
&\quad+0.5(0.5e^{-(t-\de(t))}+0.1y^2+4y^4+2y^6)\\
&=:-U_1(x,t)+0.05U_1(y,t-\de(t)).
\end{align*}
Letting $U_1(x,t)=0.5e^{-t}+0.1x^2+4x^4+2x^6$. Due to $0.9=\alpha_1(1-\bar\de)>\alpha_2=0.05$, we get condition (\ref{A3.2c}). Noting that $\varpi=1$, then condition (\ref{WF}) becomes
$$\de(t)\leq 0.08.$$

\begin{figure}[h!]
\begin{center}
\includegraphics[angle=0,height=8cm,width=9cm]{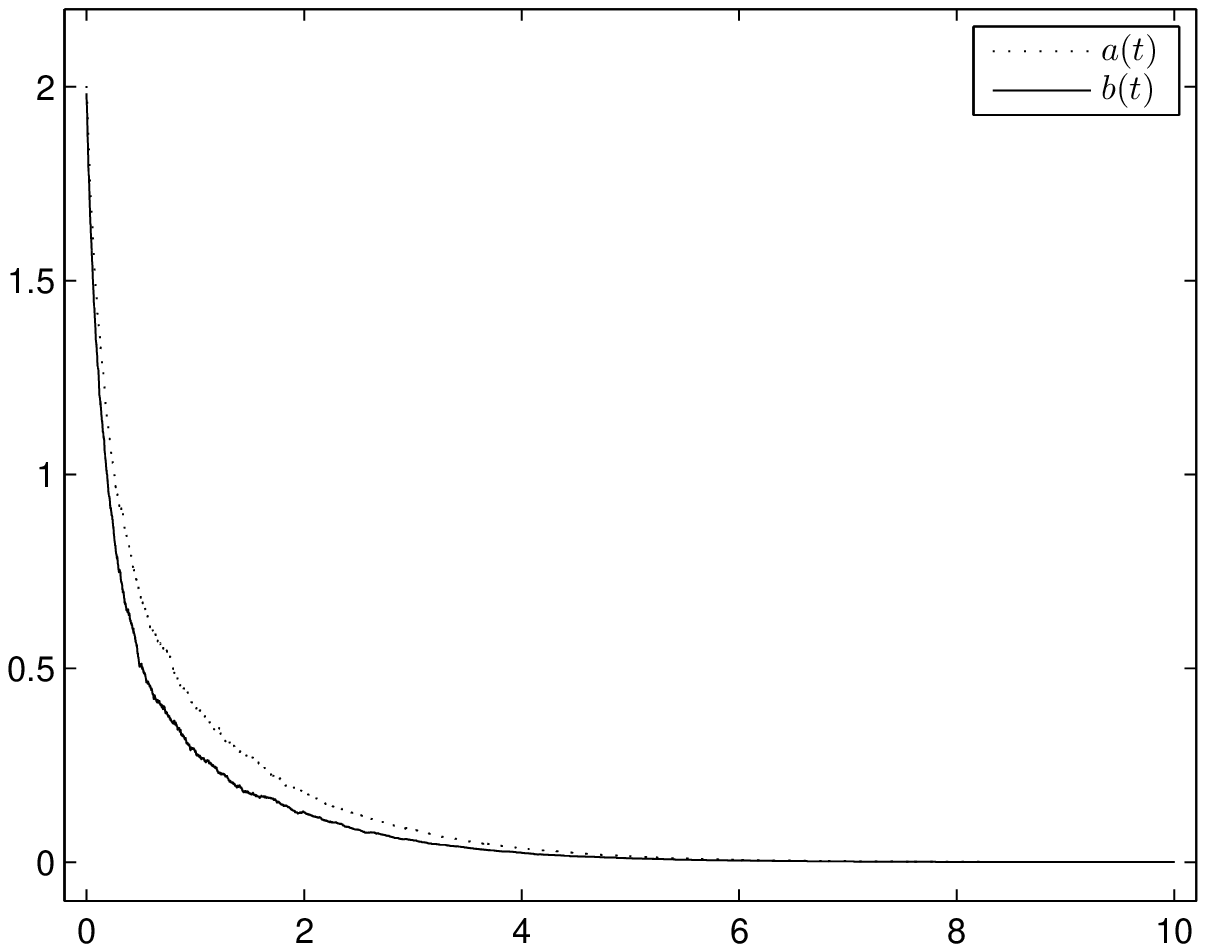}
 \end{center}
\begin{center}
\small{Figure 4.3
: The computer simulation of the upper expectation
$a(t)=\hat{\mathbb E}|X(t)|$ and the lower expectation $b(t)={\cal  E}|X(t)|$  of the solution to the $G$-SDDE  (\ref{$G$-SDDE1})  with $\de(t)=0.08$
using the Euler--Maruyama method with step size $10^{-3}$.}
\end{center}
\end{figure}

By Theorem \ref{T3.4}, we can therefore conclude that the solution of the $G$-SDDE (\ref{$G$-SDDE1}) has the properties that
\begin{align*}
&\int_0^\8 (X^2(t)+ X^4(t)+X^6(t))dt <\infty \ q.s.,\\
&{\mathcal E}\int_0^\8 (  X^2(t)+ X^4(t)+ X^6(t))dt <\infty.
\end{align*}
Moreover, as $|X(t)|^p \leq X^2(t)+ X^4(t)+ X^6(t)$ for any $p\in [2,6]$,  we have
\begin{align*}
{\mathcal E}\int_0^\8 |X(t)|^p dt <\infty.
\end{align*}
Recalling $q_1=3,\ q_2=0$ and $q=6$, we see that for $p=4$, all conditions of Theorem \ref{T3.6} are satisfied and hence we have
\begin{align*}
\lim_{t\to\infty} \hat{\mathbb{E}}|X(t)|^4   =0.
\end{align*}
We perform a computer simulation with the time-delay $\de(t)=0.08 $ for all $t\geq 0$ and the initial data $X(u)=2+\text{sin}(u)$ for $u\in [-0.08,0]$. The sample paths of the solution to the $G$-SDDE (\ref{$G$-SDDE1}) are plotted in Figure 4.3. The simulation supports our theoretical results.
 }
\end{expl}

\section{Conclusion}

  In real systems, we are often faced with two kinds of uncertainties:  probability and Knightian ones, respectively. By using Peng's sublinear expectation framework, we can characterize the systems with ambiguity. This paper gives a description of this kind of  system with delay through $G$-Brownian motion. Moreover, the boundedness and stability  of solutions to $G$-SDDEs are discussed. We first give the delay-dependent criteria of stability and boundedness of the solutions to highly nonlinear $G$-SDDEs by the method of Lyapunov functional. Then, in Section 4 we give a simulation algorithm related $G$-Brownian motion, moreover an illustrative example is explored. In fact, under classical probability, the stability of highly nonliner hybrid SDDEs has been studied by some researchers. However, our idea and framework with sublinear expectation provide a new perspective for further investigating the stability and boundedness of highly nonlinear stochastic systems driven by $G$-Brownian motion.

\section*{Acknowledgements}

The authors would also like to thank the Natural Science Foundation of China
(71571001)  for their financial support.

\end{document}